\documentclass[a4paper,leqno]{amsart}

\usepackage{amsmath}
\usepackage{amssymb}
\usepackage{amsfonts}
\usepackage[dvips]{graphicx}
\usepackage{color}
\usepackage[textheight=20cm, margin=3.7cm]{geometry}

\newtheorem{thm}{Theorem}[section]

\newtheorem{propo}[thm]{Proposition}
\newtheorem{lem}[thm]{Lemma}
\newtheorem{cor}[thm]{Corollary}

\newcommand{\R}{\mathbb{R}}
\newcommand{\C}{\mathbb{C}}
\newcommand{\Z}{\mathbb{Z}}
\newcommand{\N}{\mathbb{N}}

\newcommand{\Hh}{\mathbb{H}^2}

\newcommand{\lt}{{\mathcal L}}
\newcommand{\TO}{\mathcal{L}}
\renewcommand{\Re}{\mathrm{Re}}
\renewcommand{\Im}{\mathrm{Im}}
\DeclareMathOperator{\id}{id}
\DeclareMathOperator{\SL}{SL}
\DeclareMathOperator{\Unit}{U}

\DeclareMathOperator{\Rea}{Re}
\DeclareMathOperator{\Ima}{Im}
\DeclareMathOperator{\PSL}{PSL}

\begin{document}
\title{Fractal Weyl bounds and Hecke triangle groups}
\author[F.\@ Naud]{Fr\'ed\'eric Naud}
\address{FN: Laboratoire de Math\'ematiques d'Avignon, Universit\'e d'Avignon, 301 rue Baruch de Spinoza, 84916 Avignon Cedex, France}
\email{frederic.naud@univ-avignon.fr}

\author[A.\@ Pohl]{Anke Pohl}
\address{AP: University of Bremen, Department 3 -- Mathematics, Bibliothekstr.\@ 
5,  28359 Bremen, Germany}
\email{apohl@uni-bremen.de}

\author[L.\@ Soares]{Louis Soares}
\address{LS: Institute for Mathematics, University of Jena, Ernst-Abbe-Platz 2, 07743 Jena, Germany}
\email{louis.soares@uni-jena.de}

\begin{abstract} 
Let $ \Gamma_{w} $ be a non-cofinite Hecke triangle group with cusp width $w>2$ and let $\varrho\colon\Gamma_w\to U(V)$ be a finite-dimensional unitary representation of $\Gamma_w$. 
In this note we announce a new fractal upper bound for the Selberg zeta function of $ \Gamma_{w}$ twisted by $\varrho$. In strips parallel to the imaginary axis and bounded away from the real axis, the Selberg zeta function is bounded by $ \exp\left( C_{\varepsilon} \vert s\vert^{\delta + \varepsilon} \right) $, where $ \delta = \delta_{w} $ denotes the Hausdorff dimension of the limit set of $ \Gamma_{w}.$ 
This bound implies fractal Weyl bounds on the resonances of the Laplacian for all geometrically finite surfaces $X=\widetilde{\Gamma}\backslash \Hh$ where $\widetilde{\Gamma}$ is a finite index, torsion-free subgroup of $\Gamma_w$. 
\end{abstract}


\maketitle
\bibliographystyle{plain}

\section{Introduction and results}

Let $\Gamma$ be a finitely generated, \textit{non-cofinite}, torsion-free Fuchsian group, let $X := \Gamma\backslash\Hh$ be its associated hyperbolic surface, and let $\Delta_X$ be the hyperbolic Laplacian on $X$. 

The $L^2$-spectrum of $\Delta_X$ has been described completely by Lax--Phillips~\cite{LaxPhillips}. To be more precise, let 
$\Lambda(\Gamma)\subset \partial \Hh$ denote the limit set of $\Gamma$, and let 
\[
\delta = \delta(\Gamma) = \dim_H \Lambda(\Gamma)
\]
denote its Hausdorff dimension.  The continuous spectrum of $\Delta_X$ coincides with the half line $[1/4, +\infty)$. It does not contain embedded $L^2$-eigenvalues. The pure point spectrum is empty if $\delta\leq \frac{1}{2}$, and finite and starting at $\delta(1-\delta)$ if $\delta> \frac{1}{2}$.

This description of the spectrum implies that the resolvent 
$$R_X(s)=\left(\Delta_X-s(1-s) \right)^{-1}:L^2(X)\rightarrow L^2(X)$$
of $\Delta_X$ is well-defined and analytic on the half-plane $\{ \Rea(s)> \frac{1}{2} \}$ except at the finite set
of poles corresponding to the pure point spectrum of $\Delta_X$. The {\it resonances} of $X$ are the poles of the meromorphic continuation of 
$$R_X(s):C_0^\infty(X)\rightarrow C^{\infty}(X)$$ 
to the whole complex plane. This continuation can be deduced using the analytic Fredholm theorem after the construction of an adequate parametrix as provided by Guillop\'e--Zworski \cite{GuillopeZworski}. Throughout we denote the set of resonances of $ X $ by ${\mathcal R}_X$; each resonance is considered to be repeated according to its multiplicity. 

Understanding asymptotics for the number of resonances of $X$ is of great interest. A characterization of resonances, which proved to be helpful for such investigations, is as zeros of the Selberg zeta function $Z_X$ of $X$.

In order to be more specific, we recall that the set of prime periodic geodesics on $X$ is bijective to the set $[\Gamma]_p$ of $\Gamma$-conjugacy classes of the primitive hyperbolic elements in $\Gamma$. Let $\ell(\gamma)$ denote the length of the geodesic corresponding to $[\gamma] \in [\Gamma]_p$. For $\Re(s)>\delta$, the Selberg zeta function $Z_\Gamma$ of $X$ is given by the infinite product
\begin{equation}\label{def_szf}
Z_\Gamma(s):=\prod_{k=0}^\infty \prod_{[\gamma] \in [\Gamma]_p}\left( 1-e^{-(s+k)\ell(\gamma)}\right).
\end{equation}
For $\Re(s)\leq \delta$, the Selberg zeta function $Z_\Gamma$ is given by the meromorphic continuation of \eqref{def_szf}. The resonances are contained (including multiplicities) in the set of zeros of $Z_\Gamma$; the correspondence being exact outside of the set $\frac{1}{2}(1-\N_0)$, see \cite{BJP2}. This provides a {\it purely dynamical characterization of resonances}.

Suppose for a moment that $X=\Gamma\backslash\Hh$ has \textit{no cusps}. In this case, $\Gamma$ is a Schottky group or, equivalently, it is convex cocompact, non-cocompact, torsion-free. For these hyperbolic surfaces $X$ it was shown first by Zworski \cite{Zworski1} and then by Guillop\'{e}--Lin--Zworski \cite{Guillope_Lin_Zworski} that for all $\sigma\in\R$, as $T\rightarrow +\infty$, we have 
\begin{equation}\label{weyl_bound_schottky}
N_X(\sigma,T):=\# \{ s\in \mathcal{R}_X\ :\ \Rea(s)\geq \sigma,\ \vert \Ima(s)\vert \leq T\}=O_\sigma(T^{1+\delta}).
\end{equation}
This result is commonly referred to as a {\it fractal Weyl upper bound} by analogy with the Weyl law for eigenvalues in the compact case. These estimates were established after the pioneering work of Sj\"ostrand on semi-classical Schr\"odinger operators \cite{Sjoestrand}.

In the case that $X$ has cusps, asymptotics similar to \eqref{weyl_bound_schottky} are not yet known, mostly because of the non-compactness of the {\it trapped set} of the geodesic flow on $X$, a source of notorious difficulties. In this note, we present a fractal Weyl bound for hyperbolic surfaces arising from finite-index, torsion-free subgroups of the Hecke triangle family $\Gamma_w$.  

Hecke triangle groups are, in some sense, natural generalizations of the prominent modular group
$$
 \mathrm{PSL}_{2}(\mathbb{Z}) = \SL_2(\Z)/\{\pm \id\},
$$
which is generated by the two elements 
$$
 T:= \left (\begin{array}{cc} 1&1\\0&1 \end{array} \right ) \quad
\text{and}\quad S:= \left (\begin{array}{cc} 0&1\\-1&0 \end{array} \right ).
$$
On the hyperbolic plane $ \mathbb{H}^{2} $, these elements act by $ S(z)=-1/z $ and $ T(z) = z+1. $

The \textit{Hecke triangle group $\Gamma_w$ with cusp width $w\geq 1$} is the subgroup of $\PSL_2(\R)$ generated by $ S(z):= -1/z $ and $ T_{w}(z):= z+w $. Hecke \cite{Hecke} showed that $ \Gamma_{w} $ is Fuchsian if and only if $ w = 2 \cos\left( \pi /q \right) $ with $ q\in \mathbb{N}_{\geq 3} $, or $ w \geq 2 $. 

\begin{figure}[h]
\begin{center}
\includegraphics[scale=1.5]{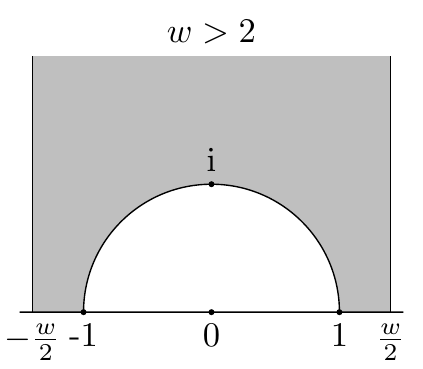}
\end{center}
\caption{Fundamental domain for $\Gamma_w$ with cusp width $w>2$.}\label{fig:funddom}
\end{figure}

In this note we focus on $\Gamma_w$ with $w>2$. For any of these non-cofinite Fuchsian groups, the associated orbifold $\Gamma_w\backslash \Hh$ is an orbifold with one cusp, one funnel and a conical singularity. The conical singularity is caused by the elliptic element $S$ in $\Gamma_w$. A fundamental domain for $\Gamma_w$ with $w>2$ is given by
$$
\mathcal{F}(w) = \left\{ z\in \mathbb{H} : \vert\mathrm{Re}(z)\vert < \tfrac{w}{2}, \ \vert z\vert > 1 \right\},
$$
see Figure~\ref{fig:funddom}.

By passing to finite index subgroups of $\Gamma_w$, one can obtain plenty of geometrically finite Fuchsian groups without elliptic elements, producing examples of Riemannian surfaces with several cusps and funnels. For example, we may consider the subgroup of index $2$ which is freely generated by the elements $T_{w}$ and $ST_{w}S$, as well as its abelian covers.

Let $\varrho:\Gamma_w\rightarrow \Unit(V)$ be a finite-dimensional unitary complex representation of $\Gamma_{w}$ with representation
space $V$. For $\Rea(s)>\delta$, the Selberg zeta function of $\Gamma_w$ twisted by $\varrho$ is given by the infinite product
$$Z_{\Gamma_w}(s,\varrho):=\prod_{k=0}^\infty \prod_{[\gamma] \in [\Gamma_w]_p}\det\left( 1_V-\varrho(\gamma)e^{-(s+k)\ell(\gamma)}\right).$$
It extends meromorphically to all of $\C$. The set of poles is contained in $\tfrac12(1-\N_0)$, and bounds on the order of poles can be given \cite{Pohl_representation, PF1}.

Our main result is the following.

\begin{thm} \label{main1}
Let $\Gamma_w$ be the Hecke triangle group with cusp width $w>2$, let $\varrho\colon\Gamma_w\to \Unit(V)$ be a finite-dimensional unitary representation of $\Gamma_w$, and let $\delta=\delta(\Gamma_w)$ denote the Hausdorff dimension of the limit set of $\Gamma_w$.  For all $\sigma_0 \in \R$, for all $\epsilon>0$, there exists $C:=C(\epsilon, \sigma_{0}, w, \varrho)$
such that for all $s=\sigma+it$ with $\sigma\geq \sigma_0$ and $\vert t \vert \geq 1$ we have
$$\log \vert Z_{\Gamma_w}(s,\varrho) \vert\leq C \vert t \vert^{\delta+\epsilon}.$$
\end{thm}

By the Venkov--Zograf factorization formula \cite{Venkov_Zograf,Venkov_book} (see also Fedosova--Pohl \cite{PF1} for the infinite-area case), the Selberg zeta function of any finite-index subgroup $\widetilde\Gamma$ of $\Gamma_w$ equals the Selberg zeta function of $\Gamma_w$ twisted with the representation of $\Gamma_w$ that is induced by the trivial one-dimensional representation of $\widetilde\Gamma$.  This allows us to deduce the following corollaries of Theorem~\ref{main1}.

\begin{cor}
\label{cor1}
Let $ w > 2 $ and let $\widetilde{\Gamma}$ be a subgroup of $\Gamma_w$ of finite index. Then 
for all $\sigma_0 \in \R$, for all $\epsilon>0$, there exists $C:=C(\epsilon,\sigma_0,\widetilde{\Gamma})$
such that for all $s=\sigma+it$ with $\sigma\geq \sigma_0$ and $\vert t \vert \geq 1$ we have
$$\log \vert Z_{\widetilde{\Gamma}}(s) \vert\leq C \vert t \vert^{\delta+\epsilon},$$
where $\delta=\delta(\Gamma_w) = \delta(\widetilde\Gamma)$.
\end{cor}

An immediate consequence of Corollary~\ref{cor1} are the following asymptotics for resonances.

\begin{cor}
\label{fractalweyl}
Let $ w > 2 $, let $\widetilde{\Gamma}$ be a finite-index subgroup of $\Gamma_w$ without elliptic elements, and set $X=\widetilde{\Gamma}
\backslash \Hh$. Then for all $\epsilon>0$, all $\sigma \in \R$, as $T\rightarrow +\infty$, we have
\begin{align*}
N_X(\sigma,T)&:=\# \{ s\in \mathcal{R}_X\ :\ \Rea(s)\geq \sigma,\ \vert \Ima(s)\vert \leq T\}=O_{\sigma,\epsilon}(T^{1+\delta+\epsilon}),
\\
M_X(\sigma,T)&:=\#\{s\in \mathcal{R}_X\ :\ \Rea(s)\geq \sigma,\ \vert \Ima(s)-T\vert \leq 1\}=O_{\sigma,\epsilon}(T^{\delta+\epsilon}).
\end{align*}
\end{cor}

Theorem \ref{main1} is the first example in the literature of a fractal Weyl bound in a situation
where the trapped set is non-compact. We expect that the techniques we developed for the family of Hecke triangle groups can be used to deal with the general geometrically finite case without major difficulties. This will be pursued in the detailed version of this work. Such bounds have direct applications: for example, one can use the techniques of Jakobson--Naud \cite{JN1} to obtain an explicit strip with infinitely many resonances. It is also reasonable to expect that such transfer operator approaches can be applied to get improved fractal Weyl bounds for resonances close to $\Rea(s)=\frac{1}{2}$, as was recently done in the convex co-compact case in \cite{Naud1,JN2, Dyatlov1}. There are also questions related to the `fractal uncertainty principle' of Dyatlov that can be understood at the level of transfer operators, see the recent paper of Dyatlov--Zworski \cite{DyatlovZworski}.

\subsubsection*{Acknowledgement} AP acknowledges support by the DFG grant PO 1483/2-1. FN is supported by Institut Universitaire de France.

\section{Outline of the proof}

Throughout let $\Gamma := \Gamma_w$ be the Hecke triangle group with cusp width $w>2$, set $ T:= T_{w} $, let $ V $ be a finite-dimensional complex vector space, endowed with a Hermitian inner product $ \langle \cdot, \cdot\rangle_{V} $ and induced norm $\|v\|_V := \sqrt{\langle v,v\rangle_V}$, and let $\varrho\colon \Gamma_w\to\Unit(V)$ be a finite-dimensional unitary representation of $\Gamma_w$. 

We start with a brief summary of the strategy that leads to a proof of Theorem~\ref{main1}. A more detailed version follows in Sections~\ref{sec:spaces}-\ref{sec:merom} below.

As we recall in Section~\ref{sec:spaces} below, the twisted Selberg zeta function $Z_{\Gamma}(s,\varrho)$ is represented by the Fredholm determinant of a twisted transfer operator family $\TO_{s,\varrho}$, $s\in\C$ (up to poles):
\[
 Z_{\Gamma}(s,\varrho) = \det(1-\TO_{s,\varrho}).
\]
Hence, it suffices to establish the estimate claimed in Theorem~\ref{main1} for $\det(1-\TO_{s,\varrho})$. To this end we may take advantage of properties of the transfer operator family itself, not just of its Fredholm determinant. There are several options for the choice of this transfer operator family. Here, we use the one given (initially only formally) by
\begin{equation}
\label{form1}
\lt_{s,\varrho}f(z):=\sum_{n\in \Z\setminus\{0\}}(\gamma'_n(z))^s\varrho^{-1}(\gamma_n) f(\gamma_n(z)), 
\end{equation}
where $\gamma_n := ST_w^n$, and $f$ belongs to a certain space of $V$-valued functions defined in a neighborhood of the limit set. For another option of the choice of the transfer operator family we refer to \cite{Pohl_representation}.

The precise choice of the function space on which $ \mathcal{L}_{s,\varrho} $ acts is crucial for our proof of Theorem~\ref{main1}. As in Guillop\'{e}--Lin--Zworski \cite{Guillope_Lin_Zworski} we use as domains of definition for $\TO_{s,\varrho}$ a certain family of vector-valued Hilbert Bergman spaces $H^2(\Omega(h);V)$ of square-integrable, holomorphic functions on a complex neighborhood $\Omega(h)$ of a certain portion $\Lambda_0$ of the limit set of $\Gamma_w$. The scaling parameter $h>0$ determines how much we `thicken' $\Lambda_0$, with smaller $h$ corresponding to less thickening. The value of $h$ is meant to be taken rather small. For our application, it is eventually optimized to essentially equal $|\Ima(s)|^{-1}$ when considering $\TO_{s,\varrho}$.

For any sufficiently small $h>0$ and any $\mathrm{Re}(s)>\tfrac12$, formula~\eqref{form1} defines $\TO_{s,\varrho}$ as an operator of trace class on $H^2(\Omega(h);V)$. Furthermore, the map $s\mapsto \TO_{s,\varrho}$ admits a meromorphic extension to all of $\C$ with poles contained in $\tfrac12(1-\N_0)$. 

For $\mathrm{Re}(s) >\tfrac12$ we can then take advantage of the standard Weyl inequality and estimate the growth of $\det(1-\TO_{s,\varrho})$ 
by the singular values of $\TO_{s,\varrho}$. To prove estimates on the singular values of $\lt_{s,\varrho}$, we adapt and extend the techniques from Bandtlow--Jenkinson \cite{Bandtlow-Jenkinson} to a vector-valued setting. 

Compared to the study of Schottky groups in \cite{Guillope_Lin_Zworski}, one of the major differences we encounter in the investigation of Hecke triangle groups, is the structure of the set $\Omega(h)$ (see Section~\ref{sec:spaces} below for its definition). As $h$ shrinks to $0$, the set $\Omega(h)$ has more and more connected components. In the case of Schottky groups, the connected components of $\Omega(h)$ have diameters uniformly of order $h$, and essentially can be considered to be Euclidean disks. In the case of Hecke triangle groups, the set $\Omega(h)$ consists of `stretched' domains with diameters ranging from $h$ to $\sqrt{h}$. To overcome the lack of uniformity in our setting, for the Hecke triangle groups, we take advantage of estimates for Patterson--Sullivan measures due to Stratmann--Urba\'{n}ski \cite{Stratmann_Urbanski} in combination with a covering argument.

For $\mathrm{Re}(s) \leq \tfrac12$, the proof of Theorem \ref{main1} is more subtle. In this region, the infinite sum \eqref{form1} no longer converges, and the transfer operator family $\TO_{s,\varrho}$ is given by its meromorphic continuation (in $s$) rather than by the infinite sum itself. This meromorphic continuation is constructive, and its main properties rely on those of the Lerch zeta function. Further details are sketched in Section~\ref{sec:merom} below. The precise knowledge of the meromorphic continuation allows us to establish sufficiently good estimates for the Fredholm determinant of $\TO_{s,\varrho}$ also in this region for $s$.

\subsection{Function spaces and transfer operators}\label{sec:spaces}

Throughout we use the Poincar\'e upper half plane model $\Hh$ for the hyperbolic plane. 
Let 
$$\Lambda_0:=\Lambda(\Gamma_w)\cap(-1,+1).$$
For any $h>0$ let $D(0,h)$ denote the open ball in $\C$ with center $0$ and radius $h$, and let
\[
\Omega(h) := \Lambda_0 + D(0,h)
\]
be the `$h$-thickening' of $\Lambda_0$ into the complex plane.  As mentioned above, the fine structure of the set $\Omega(h)$ differs 
considerably from its analogue in the case of Schottky groups. Nevertheless, the following crucial volume estimate remains true.

\begin{propo}
\label{volest}
There exists $C>0$ such that for all sufficiently small $h>0$ we have
$$\mathrm{vol}(\Omega(h))\leq Ch^{2-\delta},$$
where $ \mathrm{vol} $ denotes the Lebesgue measure on $ \mathbb{C} $, and $\delta=\delta(\Gamma_w)$ the Hausdorff dimension of the limit set $\Lambda(\Gamma_w)$ of $\Gamma_w$. 
\end{propo}

The proof of Proposition~\ref{volest} is based on an estimate, due to Stratmann--Urba\'{n}ski \cite{Stratmann_Urbanski}, of the Patterson--Sullivan measure of small intervals centered at points in $\Lambda_0$, and the Vitali covering lemma. Proposition \ref{volest} also shows that the  Minkowski dimension and the Hausdorff dimension of $\Lambda_0$ are identical.

The $V$-(vector-)valued Hilbert Bergman space on $\Omega(h)$, which we use as domain of definition for $\TO_{s,\varrho}$, is then 
$$
H^2(\Omega(h); V) :=  \left\{ \text{$f\colon \Omega(h)\to V$ holomorphic} \ \left\vert\ \Vert f\Vert_{L^{2}(\Omega(h))} < \infty \right.\right\},
$$
where 
\[
 \Vert f\Vert_{L^{2}(\Omega(h))}^{2}:=\int_{\Omega(h)} \Vert f(z)\Vert_{V}^{2} \mathrm{dvol}(z).
\]

\begin{propo}
\label{detform}
There exists $h_0>0$ such that for all $h\in (0,h_0)$ and all $s\in\C$, $\Re(s)>\frac{1}{2}$, formula~\eqref{form1} defines an operator
$$\lt_{s,\varrho}:H^2(\Omega(h);V)\rightarrow H^2(\Omega(h);V)$$
of trace class. Furthermore, 
$$Z_{\Gamma_w}(s,\varrho)=\det(1- \lt_{s,\varrho}).$$
\end{propo}

These types of identities are well-known in thermodynamical formalism. In the presence of parabolic elements, the first example was provided by Mayer \cite{Mayer_thermo}, for $\mathrm{PSL}_2(\Z)$ and twist by the trivial one-dimensional representation (or, in other words, without a twist). Work for other cofinite or non-cofinite Fuchsian groups with cusps is provided by, e.\,g., \cite{ Fried_triangle,   Morita_transfer, Mayer_Muehlenbruch_Stroemberg, Moeller_Pohl, Pohl_hecke_infinite}, an extension to non-trivial twists is shown in \cite{Pohl_representation, PF1}.

\subsection{Singular value estimates in stretched domains}\label{sec:sing}
To estimate the Fredholm determinant $\det(1-\TO_{s,\varrho})$ for $\Re(s)>\frac{1}{2}$ we  use the Weyl inequality
$$\vert \det(1-\lt_{s,\varrho})\vert \leq \prod_{k=0}^\infty (1+\mu_k(\lt_{s,\varrho})),$$
where
$$\mu_0(\lt_{s,\varrho})\geq \mu_1(\lt_{s,\varrho})\geq \ldots \geq \mu_k(\lt_{s,\varrho})\geq \ldots\geq 0$$
denotes the singular values sequence of $\lt_{s,\varrho}$, that is, the sequence of eigenvalues of the positive self-adjoint compact operator
$\sqrt{\lt_{s,\varrho}^*\lt_{s,\varrho}}$. We establish the following estimate, which is central for the proof of Theorem~\ref{main1}.

\begin{propo}
 \label{singular}
 Fix $\sigma=\Re(s)>\frac{1}{2}$, then there exist $C_1,C_2,C_3>0$ and $A>0$ such that we have for all sufficiently small $h>0$,
 $$\mu_k(\lt_{s,\varrho})\leq C_1 h^{-A}e^{C_2 \vert \Im(s)\vert h} e^{-C_3 h^\delta k}.$$
\end{propo}

For the proof of Proposition~\ref{singular} we use an extension to vector-valued situations of work by Bandtlow--Jenkinson \cite{Bandtlow-Jenkinson} for estimating singular values of general contraction operators on Hilbert Bergman spaces. A carefully designed covering argument, in which Proposition~\ref{volest} plays a crucial role, and the contraction properties of the elements $\gamma_n$ (which enable us to lift the transfer operator family to larger domains, and then take advantage of singular values estimates for canonical embeddings), then allow us to extract the explicit dependence of the bound on $h$ and $|\Ima(s)|$.

Once Proposition \ref{singular} is established, it is rather straightforward to derive Theorem \ref{main1} for $\Re(s)>\frac{1}{2}$. Indeed, setting $h=\vert \Im(s) \vert^{-1}$ we obtain for all $N\in\N$ the estimate
\begin{align*}
\log\vert \det(1- \lt_{s,\varrho})\vert
&\leq \sum_{k=0}^\infty \log(1+\mu_k(\lt_{s,\varrho}))
\\
&\leq (N+1)\log(1+C_1h^{-A})+C_1h^{-A}\sum_{k=N+1}^\infty e^{-C_3 h^\delta k}
\\
& =O(N\vert\log h \vert)+O\left(h^{-\delta-A}e^{-C_3 h^\delta N}\right).
\end{align*}
Choosing $N=C_4 h^{-\delta}\vert\log h\vert$ with $C_4>0$ sufficiently large yields for all $\epsilon>0$
$$\log\vert \det(1- \lt_{s,\varrho})\vert=O\left(h^{-\delta} \big(\log(h)\big)^2 \right)=O\left(\vert \Im(s) \vert^{\delta+\epsilon}\right),$$
and shows Theorem \ref{main1} for $ \mathrm{Re}(s) > \frac{1}{2} $. We remark that at this stage of the proof, the `$\epsilon$-loss' is caused solely by the term $h^{-A}$ in the singular values estimate. This term can be removed at the cost of a more technical proof but the improvement obtained is ruined in the next step, where the technique of meromorphic continuation will produce new polynomial losses (see Section~\ref{sec:merom} below), which we do not know how to avoid at the current state of art.

\subsection{Meromorphic continuation and Lerch zeta functions}\label{sec:merom}
The computations in Section~\ref{sec:sing} are valid for $\Re(s)>\frac{1}{2}$ only. In this region, however, the estimate from Theorem \ref{main1} is useless for asymptotics on the number of zeros of $Z_{\Gamma_w}(\cdot,\varrho)$ since it is well-known from spectral theory that in the half plane $\{ \Re(s) >\frac{1}{2} \}$, the Selberg zeta function $Z_{\Gamma_w}(\cdot,\varrho)$ has finitely many zeros only. One can even show that for the trivial one-dimensional representation ${\bf 1}_\C$, the Selberg zeta function $Z_{\Gamma_{w}}(\cdot,{\bf 1}_\C)$ has no zeros with $ \mathrm{Re}(s) > \frac{1}{2} $ except for $ s=\delta $.

To pass beyond the line $\{ \Re(s) =\frac{1}{2}\}$, we show the meromorphic continuability of $s\mapsto \lt_{s,\varrho}$ in analogy to the seminal proof by Mayer \cite{Mayer_thermo} for the case of a transfer operator family for $\SL_2(\Z)$ and using its extension to twisted transfer operators (see \cite{Pohl_representation}). 

In order to sketch this approach we recall from (\ref{form1}) that for $\Re(s)>\tfrac12$, 
$$\lt_{s,\varrho}f(z):=\sum_{n\neq 0} \frac{1}{(z+nw)^{2s}}\varrho^{-1}(\gamma_n)f\left(\frac{-1}{z+nw}\right).$$
Since $f$ is holomorphic on $\Omega(h)\ni 0$, we can use the Taylor series expansion of $f$ at $z=0$ to write (locally near $0$)
\begin{equation}\label{extend}
f(z)=f(0)+z\widetilde{f}(z),
\end{equation}
where $\widetilde{f}$ is a suitable function that is holomorphic near $0$. Using $\gamma_n$=$ST^n$ it follows that
\begin{align*}
\lt_{s,\varrho}f(z)& :=\left(\sum_{n\neq 0} \frac{1}{(z+nw)^{2s}}\varrho(T^{-n})\right) \varrho(S)f(0)
\\
& \qquad - 
\sum_{n\neq 0} \frac{1}{(z+nw)^{2s+1}}\varrho^{-1}(\gamma_n)\widetilde{f}\left(\frac{-1}{z+nw}\right)
\\
& =\left(\sum_{n\neq 0} \frac{1}{(z+nw)^{2s}}\varrho(T^{-n})\right) \varrho(S)f(0)-\lt_{s+\frac{1}{2},\varrho}\widetilde{f}(z).
\end{align*}
The term $\lt_{s+\frac{1}{2},\varrho}\widetilde f$ makes sense for $\Re(s)>0$ as well (recall that we currently suppose $\Re(s)>\tfrac12$). The first term, after choosing an appropriate basis for $V$, reads 
$$\sum_{n\neq 0} \frac{1}{(z+nw)^{2s}}\varrho(T^{-n})=\mathrm{diag}\left (  
\sum_{n\neq 0} \frac{1}{(z+nw)^{2s}}e^{2i\pi n \lambda_1},\ldots, \sum_{n\neq 0} \frac{1}{(z+nw)^{2s}}e^{2i\pi n \lambda_d} \right),$$
where $d=\mathrm{dim}(V)$ and $e^{-2i\pi\lambda_1},\ldots,e^{-2i\pi\lambda_d}$ are the eigenvalues of the unitary endomorphism $\varrho(T)$.
For $k\in\{1,\ldots,d\}$ and $\Re(s)>\frac{1}{2}$ we have
$$\sum_{n\neq 0} \frac{1}{(z+nw)^{2s}}e^{2i\pi n \lambda_k}=w^{-2s}H\left(\frac{z}{w},2s,\lambda_k\right)+w^{-2s}H\left(-\frac{z}{w}, 2s,-\lambda_k\right),$$
where 
$$H(z,s,\lambda):=\sum_{n=1}^\infty \frac{e^{2i\pi n \lambda}}{(n+z)^s}$$
is the Lerch zeta function (with a shift in the argument $z$, for convenience; compare with its definition in, e.\,g., \cite{Apostol}). It admits a meromorphic continuation in the $s$-variable to all of $\C$ with all poles contained in $1-\N_0$. 

This presentation provides a meromorphic continuation of $s\mapsto \lt_{s,\varrho}$ to the half-plane $\{ \Re(s) >0\}$, at least when the domain of $\TO_{s,\varrho}$ is restricted to functions with support in the connected component of $\Omega(h)$ that contains $0$. 

To take into account that $ \Omega(h) $ is disconnected and that we do not want to restrict the support of the functions, we use a globally defined bounded operator $\Psi_{1}\colon H^{2}(\Omega(h);V)\to H^{2}(\Omega(h);V)$ extending in a natural way the map $f\mapsto -\widetilde f$ from \eqref{extend}. We then have 
$$\lt_{s,\varrho}=\mathcal{F}_{s,\varrho,1}+\lt_{s+\frac{1}{2},\varrho}\Psi_1, $$
where $\mathcal{F}_{s,\varrho,1}$ is a rank $1$ operator essentially given by a linear combination of (finitely many) Lerch zeta functions. The procedure from above can be iterated to obtain a meromorphic continuation of $s\mapsto \TO_{s,\varrho}$ to the entire complex plane. 

\begin{propo}\label{Meromorphic}
For each $ k\in \mathbb{N} $ there exists an operator
$$ \Psi_{k} : H^{2}(\Omega(h); V)\to H^{2}(\Omega(h); V), $$
and for each $ k\in \mathbb{N} $ and $ \mathrm{Re}(s) > \frac{1}{2} $ there exists a finite-rank operator
$$ \mathcal{F}_{s,\varrho,k} : H^{2}(\Omega(h); V)\to H^{2}(\Omega(h); V) $$
such that the following holds true:
\begin{enumerate}
\item \label{Part_1} 
For all $ \mathrm{Re}(s) > \frac{1}{2} $ we have the formula
$$
\mathcal{L}_{s,\varrho} = \mathcal{F}_{s,\varrho,k} + \mathcal{L}_{s+\frac{k}{2}, \varrho} \Psi_{k}.
$$

\item \label{Part_2}
For each $ k\in \mathbb{N}$, the map  $ s\mapsto \mathcal{F}_{s,\varrho,k} $ extends to a meromorphic function on $\C$ with poles contained in $ \frac{1}{2}\left( 1- \mathbb{N}_{0}\right). $

\item \label{Part_3}
We have
$$
\Psi_{k} = \Psi_{1}^{k}, \quad  \mathcal{F}_{s,\varrho,k} = \sum_{j=0}^{k-1} \mathcal{F}_{s+\frac{j}{2}, \varrho,1} \Psi_{1}^{j}.
$$

\item \label{Part_4}
The rank of $ \mathcal{F}_{s,\varrho,k} $ is at most  $k$.
\end{enumerate}
\end{propo}

Although Proposition \ref{Meromorphic} is formulated in a rather abstract way, its proof is constructive in the sense that the operators $\Psi_{k} $ and $ \mathcal{F}_{s,\varrho, k}$ are given by explicit formulas. This enables us to control their operator norms, which is crucial for our proof of Theorem \ref{main1} for $ \mathrm{Re}(s)\leq \frac{1}{2}$.

\begin{lem}\label{Lemma_for_Psi} There exists $ h_{0}>1 $ such that for all $ h\in (0, h_{0}) $ we have
$$ \Vert \Psi_{1} \Vert\leq 5 h^{-3/2}.$$
\end{lem}  

In order to estimate the operator norm of $\mathcal{F}_{s,\varrho,k} $, we show the following growth estimate for the Lerch zeta function in strips parallel to the imaginary axis of the $ s $-plane.

\begin{propo}\label{Hurwitz_2}
Let $ \lambda\in (0,1] $ and $ 0 < r < 1 $. Write $ s $ and $ z $ in cartesian coordinates as $ s=\sigma + it $, $ z = x+iy $, and assume that $ \vert z\vert \leq r $ and $ \vert t\vert \geq 1 $. Then there exist constants $C(\sigma,\lambda,r)>0$ and $\alpha(\sigma)>0$ (with dependencies as indicated) such that  
$$ \vert H(z,s,\lambda)\vert \leq  C(\sigma,\lambda, r) \cdot \vert t\vert^{\alpha(\sigma)} \cdot e^{\frac{2}{1 - r} \vert y\vert \vert s\vert \log(1+\vert s \vert)}.  $$
\end{propo}

A key ingredient for the proof of Proposition \ref{Hurwitz_2} is a bound of the Lerch zeta function due to Katsurada \cite{Katsurada}.

Lemma \ref{Lemma_for_Psi}, Proposition \ref{Hurwitz_2}, and Proposition \ref{Meromorphic}\eqref{Part_3} imply that for all sufficiently small $ h > 0 $ we have 
\begin{equation}\label{Psi_lemma}
\Vert \Psi_{k}\Vert \leq 5^{k}h^{-3k/2}
\end{equation}
and 
\begin{equation}\label{F_lemma}
\Vert \mathcal{F}_{s,\varrho,k}\Vert\leq C h^{-C}\vert \mathrm{Im}(s)\vert^{C} e^{ C h ( \vert s\vert +\frac{k}{2}) \log(1+\vert s \vert +\frac{k}{2})}
\end{equation}
for some constant $ C = C(\Re(s), k, w,\varrho)>0$.

To estimate the Fredholm determinant of $\TO_{s,\varrho}$ in the region of meromorphic continuation we do not use the Weyl inequality directly as in Section~\ref{sec:sing}, but the following refined version of it. The proof of it is fairly simple and follows from an inequality among determinants by Seiler--Simon \cite{SS}.

\begin{lem}\label{Determinant_Inequality}
Let $ \mathcal{H} $ be a separable Hilbert space, $ \mathcal{F} : \mathcal{H}\to \mathcal{H} $ a finite-rank operator, and $ \mathcal{T} : \mathcal{H}\to \mathcal{H} $ an arbitrary trace class operator. Then
$$ \log \vert \det(1+\mathcal{F}+\mathcal{T})\vert \leq \mathrm{rank}(\mathcal{F}) \log(1+\Vert \mathcal{F}\Vert)+\sum_{m=1}^{\infty} \log(1+\mu_{m}(\mathcal{T})). $$
\end{lem}

To complete the proof of Theorem~\ref{main1}, we fix $\sigma=\Re(s)$ and choose $k\in \mathbb{N}$ sufficiently large such that $\sigma+\frac{k}{2}>\frac{1}{2}$. Applying Lemma \ref{Determinant_Inequality} gives
$$\log \vert Z_{\Gamma_w}(s,\varrho) \vert = \log \vert\det (1-\mathcal{F}_{s,\varrho,k}-\mathcal{L}_{s+\frac{k}{2}, \varrho} \Psi_{k})\vert $$
$$\leq k\log(1+\Vert \mathcal{F}_{s,\varrho,k}\Vert)+
\sum_{m=0}^\infty \log\left(1+\mu_m(\mathcal{L}_{s+\frac{k}{2}, \varrho})\Vert \Psi_k \Vert \right).$$
Invoking the norm estimates \eqref{Psi_lemma} and \eqref{F_lemma} as well as the singular values estimate from Proposition \ref{singular}, we obtain 
$$\log \vert Z_{\Gamma_w}(s,\varrho) \vert\leq O(\vert \log h\vert)+\sum_{m=0}^\infty \log\left(1+Ch^{-A'}e^{-Ch^\delta m} \right),$$
for some $C,A'>0$. Then one continues as in Section~\ref{sec:sing}.

We remark that a more careful choice of $h$ yields the more accurate bound
\begin{equation}\label{last_bound}
\log \vert Z_{\Gamma_w}(s,\varrho) \vert\leq O\left( \vert\Im(s)\vert^{\delta} (\log \vert \Im(s)\vert )^{2-\delta}\right).
\end{equation}
However, our methods do not allow the `log-loss' to be removed entirely, and we do not know if the upper bound \eqref{last_bound} is optimal.

\subsection{Conclusions}
We have announced and sketched the proof of fractal Weyl bounds for resonances of geometrically finite hyperbolic surfaces arising from torsion-free, finite index subgroups of Hecke triangle groups. In the detailed version, which will be published elsewhere, we will generalize these estimates to a much wider class of geometrically finite surfaces for which similar transfer operators families can be constructed. 
The ideas remain the same, but the level of generality adds technical and notational complexity to the proof.

\end{document}